\documentstyle[12pt]{article}
\titlepage

\title{When is a Connection a Levi-Civita Connection?}
\author{  Richard Atkins \\
        richard.atkins@twu.ca \\
		Department of Mathematics\\ 
		Trinity Western University \\
		7600 Glover Road \\
		Langley, BC, V2Y 1Y1 Canada}
\date{}
\newtheorem{fact}{Fact}

\newtheorem{lemma}[fact]{Lemma}
\newtheorem{proposition}[fact]{Proposition}
\newtheorem{theorem}[fact]{Theorem}
\newtheorem{corollary}[fact]{Corollary}
\begin{document}
\maketitle
\begin{abstract}
We consider the more general question as to when a connection is a metric connection. 
There are two aspects to this investigation: first, 
the determination of the integrability conditions that ensure the existence 
of a local parallel metric in the neighbourhood of a given point
and second, the characterization of the topological obstruction to a 
globally defined parallel metric. \\
\begin{center} MSC: 53C05
\end{center}
\newpage
By means of a derived flag it may be 
ascertained when a connection on a vector bundle possesses, locally, a non-trivial 
parallel section; in particular, for the bundle of symmetric, 
covariant two-tensors on a manifold this solves
the problem as to when a connection is locally metric. For connections satisfying a certain
regularity requirement, the existence of a global parallel metric may be formulated in 
the setting of the cohomology of the constant sheaf of sections in the general linear group. 
When the general linear group associated to the connection is one-dimensional the 
obstruction may be defined in terms of de Rham cohomology. Lastly, the difference between
analytic and smooth connections is examined: whereas analytic connections are regular,
the case of smooth connections is shown to possess a pathology in that a locally metric 
connection need not be globally metric, even on a contractible manifold.
\end{abstract}
\newpage
\section{Introduction}
A connection on a manifold is a type of differentiation that acts on vector fields,
differential forms and tensor products of these objects. 
Its importance lies in the fact that given a piecewise continuous curve connecting two 
points on the manifold, the connection defines a linear isomorphism between the respective 
tangent spaces at these points.  

Another fundamental concept in the study of differential geometry is that of a metric,
which provides a measure of distance on the manifold. It is well known that a metric 
uniquely determines a Levi-Civita connection: a symmetric connection for which the metric 
is parallel. Not all connections are derived from a metric and so it is natural to ask,  
for a symmetric connection, if there exists a parallel metric, that is, 
whether the connection is a Levi-Civita connection. More generally, a connection
on a manifold $M$, symmetric or not, is said to be {\it metric} if there exists a parallel 
metric defined on $M$. Henceforth all manifolds are assumed to be connected since 
the question of the existence of a parallel metric on a manifold reduces to the 
corresponding problem on its connected components. 

One may also be concerned with the possibility of local parallel metrics on the manifold.
Specifically, a connection $\nabla$ on $M$ is {\it locally metric at $x\in M$} if there exists 
a neighbourhood $U$ of $x$ and a metric defined on $U$, which is parallel with respect to the 
connection restricted to $U$. $\nabla$ is {\it locally metric} if it is locally metric
at each point $x\in M$.  
The investigation into whether a connection is metric comprises two aspects: first, 
the determination of the integrability conditions that ensure the existence of a local 
parallel metric in the neighbourhood of a given point and second, the characterization of 
the topological obstruction to a globally defined parallel metric. 

In the next section it is determined when a connection on an arbitrary vector bundle
possesses, locally, a non-trivial parallel section. 
This is accomplished by means of calculating
a derived flag of subsets of the bundle. If the terminal subset $\widetilde{W}$ admits, 
locally, the structure of a non-zero vector bundle then local non-trivial 
parallel sections exist. By considering the vector bundle of symmetric, covariant two-tensors on a 
manifold the question as to whether the connection is locally metric is answered.

Section 3 considers the problem of the existence of a global parallel 
metric when $\widetilde{W}$ is a vector bundle of rank one over $M$. 
If $\nabla$ is a locally metric connection then the existence of a global 
parallel metric may be expressed in terms of de Rham cohomology.

In the two sections following, the connection is required to be {\it regular} 
in the sense that $\widetilde{W}$ is a rank $n$ vector bundle. As such it is 
a flat vector bundle with respect to the connection. It is known that the group of
isomorphism classes of
flat line bundles over a manifold is isomorphic to the first cohomology group of the 
constant sheaf of sections in the multiplicative group of non-zero complex numbers
(cf. \cite{cc}, 2.2.11. Theorem (3), pg. 76). In Section 4, this is generalized to the 
classification of flat vector bundles of arbitrary rank. This leads to the study of the
cohomology $H^{1}(M,{\cal G}_{n})$ of the constant sheaf of sections in the 
general linear group $GL(n,{\Re})$. In Section 5, the question as to whether a 
connection is metric is answered for the case of regular connections. 

In the final section, the difference between analytic and smooth connections is explored.
Analytic connections are always regular and so admit the cohomological description
developed in the previous two sections. In particular, an analytic locally metric connection
on an analytic simply-connected manifold is metric. This is not the case with 
smooth connections in general. We close with an example of a smooth
locally metric connection on
the plane, which is not globally metric. Thus even with contractible spaces the local parallel
metrics might not piece together to yield a global parallel metric.

Similar work has been investigated in \cite{aa} in the case of surfaces and \cite{dd} for
four-dimensional manifolds.

\section{When is a connection locally metric?}

We first consider the more general question as to when a connection on a vector bundle
admits, locally, a non-trivial parallel section.

Let $\pi: W\rightarrow M$ be a vector bundle and $W'$ a subset of $W$ with the 
following two properties: \\
P1: For each $x\in M$,  $W_{x}\cap W'$ is a linear subspace of $W_{x}:=\pi^{-1}(x)$. \\
P2: For each $ w\in W'$ there exists an open neighbourhood $U$ of $\pi(w)$ in $M$
and a smooth local section $X: U\subseteq M\rightarrow W'\subseteq W$ such that
$w = X(\pi(w))$. 

Let \[ \nabla :{\cal A}^{0}(W) \rightarrow {\cal A}^{1}(W) \]
be a connection on $W$, where ${\cal A}^{n}(W)$ denotes the space of local sections 
$U\subseteq M \rightarrow W \otimes \Lambda^{n}M$. Define a map 
\[ \widetilde{\alpha}: {\cal A}^{0}(W) \rightarrow {\cal A}^{1}(W/W') \] 
by \[ \widetilde{\alpha} : = \phi \circ \nabla, \]
where $W/W'$ is the quotient of $W$ and $W'$ taken fibrewise and
$\phi:W\otimes T^{*}M \rightarrow (W/W')\otimes T^{*}M$ denotes the natural projection.
For any local section $X:U\subseteq M \rightarrow W'\subseteq W$ and 
differentiable function $f:U\rightarrow \Re$
we have $\widetilde{\alpha}(fX) = f\widetilde{\alpha}(X)$. Thus, there corresponds to 
$\widetilde{\alpha}$ a map 
\[ \alpha_{W'} : W' \rightarrow (W/W')\otimes T^{*}M \]
acting linearly on each fibre of $W'$.
$\alpha_{W'}$ is the {\it second fundamental 1-form} of $W'$. 

Let $V$ be any subset of $W$ satisfying P1. Define ${\cal S}(V)$ to be the subset of
$V$ consisting of all elements $v$ for which there exists a smooth local section
$X:U\subseteq M \rightarrow V\subseteq W$ such that $v = X(\pi(v))$. Then 
${\cal S}(V)$ satisfies both P1 and P2.

We seek to construct the maximal flat subset $\widetilde{W}$,  of $W$. 
$\widetilde{W}$ may be obtained as follows. Set
\[ \begin{array}{lll}
   V^{(0)} & := & \{ w \in W \hspace{.03in} | \hspace{.03in} R(,)(w) = 0 \} \\
   W^{(i)}   & := & {\cal S}(V^{(i)}) \\
   V^{(i+1)} & := & ker \hspace{.03in} \alpha_{W^{(i)}} 
  \end{array} \]
where $R:TM\otimes TM\otimes W\rightarrow W$ 
denotes the curvature tensor of $\nabla$. This gives a sequence
\[ W \supseteq W^{(0)} \supseteq W^{(1)} \supseteq \cdots \supseteq W^{(k)} 
\supseteq \cdots \]
of subsets of $W$. Note that $W^{(i)}$ is not necessarily a vector bundle over $M$ since
the dimension of the fibres may vary from point to point.
For some $k \in { N}$, $W^{(l)}= W^{(k)}$ for all $l \geq k$. 
Define $\widetilde{W} = W^{(k)}$, with projection $\tilde{\pi}:\widetilde{W} \rightarrow M$. 

In order to extract information from $\widetilde{W}$ we need some concept of regularity. 
Accordingly, we say that the connection $\nabla$ is {\it regular at $x\in M$} if there exists 
a neighbourhood $U$ of $x$ such that $\tilde{\pi}^{-1}(U) \subseteq \widetilde{W}$
is a vector bundle over $U$. $\nabla$ is $regular$ if $\widetilde{W}$ is a vector bundle over
$M$. The dimension of the fibres of $\widetilde{W}$, for regular $\nabla$,
shall be denoted $rank \hspace{0.03in} \widetilde{W}$. 

\begin{theorem} \label{theorem:flag}
Let $\nabla$ be a connection on the vector bundle $\pi:W\rightarrow M$. \\
(i) If $X:U \subseteq M \rightarrow W$ is a local parallel section then the image of $X$ lies in
$\widetilde{W}$. \\
(ii) Suppose that $\nabla$ is regular at $x\in M$. Then for every $w \in \widetilde{W}_{x}$
there exists a local parallel section 
$X:U\subseteq M \rightarrow \widetilde{W}$ with $X(x) = w$.
\end{theorem}
{\bf Proof:}\\
(i) follows directly from the definition of $\widetilde{W}$. \\
(ii) Suppose that $\nabla$ is regular at $x\in M$ and let $w \in \widetilde{W}_{x}$. 
By regularity, there exists a neighbourhood $U_{1}$ of $x$
and a frame $(X_{1},...,X_{n})$ of $\tilde{\pi}^{-1}(U_{1})\subseteq \widetilde{W}$.
By choosing a possibly smaller
neighbourhood $U_{2} \subseteq U_{1}$ of $x$ we can extend $(X_{1},...,X_{n})$ to a  
frame ${\cal X} := (X_{1},...,X_{n}, ...,X_{N})$ of $\pi^{-1}(U_{2})\subseteq W$. 
Let $\omega = \omega^{i}_{j}$ denote
the connection form of $\nabla$ with respect to  
${\cal X}$: $\nabla_{X}X_{j}=\sum_{i=1}^{N}X_{i} \omega^{i}_{j}(X)$.
Since $\widetilde{W}$ has zero second fundamental 1-form,
\[ \omega = \left( \begin{tabular}{c|c} 
$\phi$ & $*$ \\ \hline
$0$ & $*$  \end{tabular} \right) \]
where $\phi$ is an $n \times n$ matrix of 1-forms. The curvature
form $\Omega = \Omega^{i}_{j}$ of $\nabla$ with respect to ${\cal X}$ is

\[ \Omega = d\omega + \omega\wedge\omega = \left( \begin{tabular}{c|c}
$d\phi+\phi\wedge\phi$ & $*$ \\ \hline
$0$ & $*$ \end{tabular} \right) \]
Since the curvature tensor $R$ is identically zero, when restricted to
$TM\otimes TM\otimes \widetilde{W}$, it follows that
\[ d\phi+\phi\wedge\phi = 0\]
 
Therefore, by the Frobenius Theorem, 
there exists an $n \times n$ matrix of functions $A= A^{i}_{j}$
defined in a 
neighbourhood $U\subseteq U_{2}$ of $x$ such that 
$dA = -\phi \wedge A$ and $A(x) = I_{n \times n}$, 
the $n \times n$ identity matrix (cf. \cite{ee}, chp. 7, 2. Proposition 1., 
pg. 290). Let $c^{j}$, $1\leq j \leq n$, be real scalars satisfying 
$w = \sum_{j=1}^{n}X_{j}(x) c^{j}$. Define functions $f^{i}$ on $U$  
by 
\[ f^{i} = \left\{ \begin{array}{ll}
\sum_{j=1}^{n} A^{i}_{j}c^{j} & \hspace{0.5in} 1 \leq i\leq n \\
0 & \hspace{0.5in} n+1 \leq i \leq N. \end{array} \right. \]
Let $X:U\rightarrow \widetilde{W}$ be the local section of $\widetilde{W}$ defined  by
$X:= \sum_{i=1}^{N} X_{i}f^{i}$. Since $df+\omega\cdot f=0$, $X$ is  parallel.
Moreover, $X(x)=w$. \\
{\bf q.e.d.}

\begin{corollary} \label{corollary:flat}
Let $\nabla$ be a regular connection on the vector bundle $\pi:W\rightarrow M$.
Then $(\widetilde{W}, \nabla)$ is a flat vector bundle over $M$.
\end{corollary}

\begin{corollary} \label{corollary:number}
Let $\nabla$ be a connection on the vector bundle $\pi:W\rightarrow M$, regular at
$x\in M$. Then there are  $dim \hspace{0.03in} \widetilde{W}_{x}$ independent local 
parallel sections in a neighbourhood of $x\in M$. 
\end{corollary}

Now let $W$ be the vector bundle over $M$ consisting of the symmetric elements of 
$T^{*}M\otimes T^{*}M$. The connections that are locally metric are characterized as follows.

\begin{corollary} \label{corollary:local}
Let $\nabla$ be a connection on $M$, regular at $x\in M$. Then
$\nabla$ is locally metric at $x$ if and only if  $\widetilde{W}_{x}$ 
contains a positive-definite bilinear form.
\end{corollary}

The Frobenius Theorem, upon which Theorem \ref{theorem:flag} is based, has both smooth
and real analytic versions. Therefore Theorem \ref{theorem:flag} and Corollaries
\ref{corollary:flat}, \ref{corollary:number} and \ref{corollary:local} may be interpreted 
either in the smooth or analytic contexts. Henceforth, however, manifolds and 
connections shall be assumed to be smooth unless explicitly stated to be analytic. \\

\hspace{-0.3in} {\bf Example 1} Consider the symmetric connection $\nabla$ on the 2-sphere, 
$M=S^{2}$, defined as follows:
$\Gamma^{\theta}_{\phi \phi}=-sin{\theta}cos{\theta}$, $\Gamma^{\phi}_{\theta \phi}=
\Gamma^{\phi}_{\phi \theta}=cot{\theta}$ and all other Christoffel symbols are zero. 
Here $\theta$ and
$\phi$ are the polar and azimuthal angles on $S^{2}$, respectively. Let 
\[ \begin{array}{lll}
   X_{1} & = & d\theta \otimes d\theta \\
   X_{2} & = & d\phi \otimes d\phi \\
   X_{3} & = & d\theta \otimes d\phi + d\phi \otimes d\theta 
  \end{array} \]
be a basis of $W$, the symmetric elements of $T^{*}M\otimes T^{*}M$.
The curvature terms 
$R_{\theta \phi}=\nabla_\frac{\partial}{\partial \theta}
\nabla_\frac{\partial}{\partial \phi}-\nabla_\frac{\partial}{\partial  \phi}
\nabla_\frac{\partial}{\partial \theta}$ are
\[ \begin{array}{lll}
   R_{\theta \phi}(X_{1}) & = & -(sin^{2}{\theta})X_{3}  \\
   R_{\theta \phi}(X_{2}) & = & X_{3}  \\ 
   R_{\theta \phi}(X_{3}) & = & 2X_{1}-2(sin^{2}{\theta})X_{2}  
     \end{array} \]
This gives $W^{(0)} = span(X_{1} +(sin^{2}{\theta}) X_{2})$. 
Non-zero local sections of $W^{(0)}$ are of the form 
$X=f(X_{1}+(sin^{2}{\theta}) X_{2})$ where $f$ is a smooth non-vanishing function defined 
on an open subset of 
$S^{2}$. The covariant derivative of $X$ is $\nabla X=X \otimes dlog|f|$ and so 
$W^{(1)} = W^{(0)}$. Thus $\widetilde{W} = W^{(0)}$. 
Since $\widetilde{W}$ is a rank one vector bundle over $S^{2}$  it follows that $\nabla$
is a locally metric connection.

\section{The global problem: $rank \hspace{0.03in} \widetilde{W}  = 1$}

Consider a locally metric connection $\nabla$ on a manifold $M$ such that 
$\widetilde{W}$ is a rank one vector bundle over $M$. In such a case 
the existence of a global parallel metric on $M$ can be expressed in terms 
of de Rham cohomology.

Let $\widetilde{W}^{+}$ denote the subset of $\widetilde{W}$ consisting of the 
positive-definite bilinear forms. Since $\nabla$ is locally metric, $\widetilde{W}^{+}$
is a rank one fibre bundle over $M$.
Let $s:M\rightarrow \widetilde{W}^{+}$ be any section 
of $\widetilde{W}^{+}$. $s$ is a metric on $M$, which is not necessarily parallel. Consider
an open cover $\{U_{\alpha}: \alpha \in A \}$ of $M$ for which there exists a local parallel 
metric $h_{\alpha}$ on $U_{\alpha}$, for each $\alpha \in A$. Then $h_{\alpha}$ is a 
local section: $h_{\alpha}:U_{\alpha}\rightarrow \widetilde{W}^{+}$. 
Since $\widetilde{W}^{+}$ is a rank one fibre bundle, $s$ restricted to $U_{\alpha}$, 
denoted $s_{\alpha}$, 
must be of the form $s_{\alpha}=f_{\alpha}h_{\alpha}$ for some smooth, positive function 
$f_{\alpha}$ on $U_{\alpha}$. This allows us to 
define the 1-form $\Phi_\alpha$ on $U_{\alpha}$ by 
the covariant derivative of $s_{\alpha}$:
$\nabla s_{\alpha} = s_{\alpha} \otimes \Phi_{\alpha}$. 
Since $s_{\alpha}=s_{\beta}$ on $U_{\alpha} \cap U_{\beta}$ it follows that $\Phi_{\alpha} =
\Phi_{\beta}$ on $U_{\alpha} \cap U_{\beta}$ as well. Therefore there exists a unique 1-form
$\Phi$ on $M$ defined by $\nabla s = s \otimes \Phi$. $\Phi$ restricted to $U_{\alpha}$ is just
$\Phi_{\alpha} = dlogf_{\alpha}$, which is exact. Therefore $\Phi$ is a closed 1-form on $M$ and thus 
defines a cohomology class $[\Phi]$ in $H^{1}_{de R}(M)$. 

It is easily seen that $[\Phi]$ does not depend on the choice of section $s$.
Indeed, let $s':M\rightarrow \widetilde{W}^{+}$ be another section. Then $s'=fs$ for some
positive function $f$ on $M$. $\nabla s'=s'\otimes\Phi'$, where $\Phi'= \Phi+dlogf$.
Thus $[\Phi']=[\Phi]$. We may now state when $\nabla$ is a metric connection.

\begin{theorem} \label{theorem:dim1}
Let $\nabla$ be a regular, locally metric connection on $M$ with  
$rank \hspace{0.03in} \widetilde{W}  = 1$. Then
$\nabla$ is a metric connection on $M$ if and only if $[\Phi]=0$ in $H^{1}_{de R}(M)$.
\end{theorem}
{\bf Proof:}\\ 
$\Longrightarrow$ Suppose that $\nabla$ is a metric connection on $M$. That is, there 
exists a parallel metric $s$, which must be a section $s:M\rightarrow \widetilde{W}^{+}$. 
Then $\nabla s=0$ defines $\Phi=0$.

\hspace{-0.3in} $\Longleftarrow$ Suppose that $\Phi$ is exact. We may write $\Phi = df$ for 
some function $f$ on $M$. Define the metric $h:=exp(-f)s$. $h$ is parallel and so $\nabla$
is metric. \\ 
{\bf q.e.d.} \\

\hspace{-0.3in} {\bf Example 1 continued}  A section of $\widetilde{W}^{+}$ is given by 
$s:= X_{1}+(sin^{2}{\theta})X_{2}$. This defines the 1-form $\Phi=0$ which is trivially exact. 
Therefore the connection $\nabla$ is, in fact, a Levi-Civita connection on $S^{2}$; 
it is the Levi-Civita connection of the metric induced by the standard embedding of $S^{2}$ 
into $\Re^{3}$.

\section{Classification of flat vector bundles}

We have seen (Corollary \ref{corollary:flat}) that for a regular connection $\nabla$ on
$M$, $(\widetilde{W}, \nabla)$ is a flat vector bundle. In pursuit, 
therefore, of the topological obstruction to the existence of a global parallel metric
it will be convenient to classify the flat vector bundles over $M$, 
for which the language of sheaf cohomology will prove useful. 
This shall be the focus of the present section. 

Let ${\cal U} = \{U_{\alpha}: \alpha \in A \}$ 
be an open cover of $M$.  Let ${\cal G}_{n}$ denote the constant sheaf of sections 
in the group $GL(n; \Re)$. Consider the \v{C}ech complex 
$C^{*}({\cal U},{\cal G}_{n}):=(C^{0}({\cal U},{\cal G}_{n}), C^{1}({\cal U},{\cal G}_{n}),
C^{2}({\cal U},{\cal G}_{n}))$. Recall that the coboundary operators for non-abelian cohomology
$\delta^{i} : C^{i}({\cal U},{\cal G}_{n}) \rightarrow C^{i+1}({\cal U},{\cal G}_{n})$ 
are defined by
\[ \begin{array}{lll}
   (\delta^{0} a)_{\alpha \beta} & := & a_{\alpha}^{-1} a_{\beta}  \\
   (\delta^{1} a)_{\alpha \beta \gamma}  & := & 
     a_{\alpha \gamma}^{-1} a_{\alpha \beta} a_{\beta \gamma}   
   \end{array} \]
An element $a \in C^{1}({\cal U},{\cal G}_{n})$ is a {\it 1-cocycle}
if $\delta^{1} a = 1$, where
$1_{\alpha \beta \gamma}$ is the $n\times n$ identity matrix. 
$a$ is a {\it 1-coboundary} if there
exists an element $b \in  C^{0}({\cal U},{\cal G}_{n})$ such that $a=\delta^{0} b$.
Define an equivalence relation $\sim$ on the set of 1-cocyles by
$a\sim a'$ iff there exists an element $b$ in $C^{0}({\cal U},{\cal G}_{n})$ such that 
$a'_{\alpha \beta} =b_{\alpha}^{-1} a_{\alpha \beta} b_{\beta}$. Equivalent elements are
said to be {\it cohomologous} and the set of elements cohomologous to $a$ is denoted $[a]$. 
The set of equivalence classes defines the 
{\it first cohomology set $\check{H}^{1}({\cal U}, {\cal G}_{n})$ of the complex 
$C^{*}({\cal U},{\cal G}_{n})$ }.
Observe that it is not a group, in general, but a pointed set with distinguished element  
$1=1_{n}$. 
The {\it first cohomology set} is defined to be the direct limit 
\[ H^{1}(M, {\cal G}_{n}) := \lim_{\rightarrow} \check{H}^{1}({\cal U}, {\cal G}_{n}) \]

Denote by ${\cal E}_{n}$ the set of isomorphism classes of
pairs $(E,\nabla)$, where $E\rightarrow M$ is a rank $n$ real vector bundle 
over $M$ and $\nabla$ is a flat connection on $E$.

\begin{proposition} \label{proposition:flat} 
\hspace{0.1in} \\
(i) Each isomorphism class $\xi \in {\cal E}_{n}$ determines an element $\Theta(\xi) \in
 H^{1}(M, {\cal G}_{n})$. \\
(ii) Each element $\theta\in H^{1}(M, {\cal G}_{n})$ determines an isomorphism class
$\Psi(\theta) \in {\cal E}_{n}$. \\
(iii) $\Theta(\Psi(\theta)) = \theta$. \\
(iv) $\Psi(\Theta(\xi)) = \xi$.\\
\end{proposition}
{\bf Proof:}\\ 
(i) Let $(E,\nabla)$ be a representative of $\xi \in {\cal E}_{n}$.
There exists an open cover ${\cal U}=\{U_{\alpha}: \alpha \in A\}$ of $M$ with associated
bases of local parallel sections $h_{\alpha i}: U_{\alpha}\rightarrow E$, $1\leq i \leq n$. 
On non-empty intersections of the open sets, $U_{\alpha \beta}:= U_{\alpha}\cap U_{\beta}$, 
the bases  on $U_{\alpha}$ and $U_{\beta}$ are related by : 
\[ h_{\alpha i |_{U_{\alpha \beta}}} = \sum_{j=1}^{n}(t_{\alpha \beta})_{i j}
   h_{\beta j |_{U_{\alpha \beta}}}, \]
for some $t\in C^{1}({\cal U},{\cal G}_{n})$, for $1 \leq i \leq n$. 
The 1-cochain $t \in C^{1}({\cal U},{\cal G}_{n})$ satisfies 
\[ (\delta^{1} t)_{\alpha \beta \gamma} := 
           t_{\alpha \gamma}^{-1} t_{\alpha \beta} t_{\beta \gamma}  = 1_{\alpha
            \beta \gamma} \]
Therefore $t$ is a 1-cocycle and thus defines an element 
$\tau:= [t]$ in $H^{1}(M, {\cal G}_{n})$. 

Let ${\cal U'}=\{U'_{\alpha}: \alpha \in A'\}$ be another open cover of $M$ with 
associated bases 
of parallel sections $h'_{\alpha i}: U_{\alpha}\rightarrow E$, $1\leq i \leq n$. The 
element $t' \in C^{1}({\cal U'},{\cal G}_{n})$ giving the change of basis over 
non-empty intersections
$U'_{\alpha \beta}:= U'_{\alpha}\cap U'_{\beta}$ also defines an element $\tau':=[t']$ in
$H^{1}(M, {\cal G}_{n})$. It must be shown that $\tau' = \tau$. Let 
${\cal V}=\{V_{\alpha}: \alpha \in B \}$  be a common refinement of ${\cal U}$ and ${\cal U'}$: 
$V_{\alpha} \subset U_{\rho(\alpha)}$ and  $V_{\alpha} \subset U'_{\rho'(\alpha)}$, where
$\rho:B\rightarrow A$ and $\rho':B\rightarrow A'$. Let $\hat{t}$ (resp. $\hat{t}'$)
denote the image of $t$ (resp. $t'$) in $C^{1}({\cal V},{\cal G}_{n})$:
$\hat{t}_{\alpha \beta} = t_{\rho(\alpha) \rho(\beta)|_{V_{\alpha \beta}}}$ and 
$\hat{t}'_{\alpha \beta} = t'_{\rho'(\alpha) \rho'(\beta)|_{V_{\alpha \beta}}}$,
on non-empty intersections $V_{\alpha \beta} := V_{\alpha} \cap V_{\beta}$.
There exists $a\in C^{0}({\cal V},{\cal G}_{n})$ giving the change of basis:
$h_{\rho(\alpha)i|_{V_{\alpha}}} = 
\sum_{j=1}^{n} (a_{\alpha})_{ij}h'_{\rho'(\alpha)j|_{V_{\alpha}}}$. Thus
\begin{eqnarray*}
h'_{\rho'(\alpha)i|_{V_{\alpha \beta}}}
& = \hspace*{.1in} & 
\sum_{j=1}^{n}(a_{\alpha})^{-1}_{ij}h_{\rho(\alpha) j|_{V_{\alpha \beta}}} \\
& = \hspace*{.1in} & 
\sum_{j,k=1}^{n}(a_{\alpha})^{-1}_{ij} (\hat{t}_{\alpha \beta})_{jk}
h_{\rho(\beta) k|_{V_{\alpha \beta}}}   \\
& = \hspace*{.1in} & 
\sum_{j,k,l=1}^{n}(a_{\alpha})^{-1}_{ij} (\hat{t}_{\alpha \beta})_{jk}
(a_{\beta})_{kl} h'_{\rho'(\beta) l|_{V_{\alpha \beta}}}  
\end{eqnarray*}
Since 
\[ h'_{\rho'(\alpha)i|_{V_{\alpha \beta}}} = 
\sum_{l=1}^{n}(\hat{t}'_{\alpha \beta})_{il}
h'_{\rho'(\beta) l|_{V_{\alpha \beta}}} \]
we have $\hat{t}'_{\alpha \beta} = a_{\alpha}^{-1}\hat{t}_{\alpha \beta}a_{\beta}$. 
Hence $\hat{t}'$ is cohomologous
to $\hat{t}$ in $\check{H}^{1}({\cal V}, {\cal G}_{n})$ and so $\tau' = \tau$
in $H^{1}(M, {\cal G}_{n})$. Therefore $\Theta$ is well-defined.\\

\hspace{-0.25in}(ii) Let $\theta \in H^{1}(M, {\cal G}_{n})$. $\theta$ is represented by a 
1-cocycle $t \in C^{1}({\cal U},{\cal G}_{n})$, for some cover 
${\cal U}= \{ U_{\alpha}: \alpha \in A \}$
of $M$. The elements $t_{\alpha \beta}$ are the locally constant transition 
functions of an associated rank $n$ real vector bundle $\pi: E\rightarrow M$. 
On each $\pi^{-1}(U_{\alpha})$ there is a canonical basis 
$e_{\alpha i}$, $ 1\leq i \leq n$, of vector fields on $U_{\alpha}$ and on 
$U_{\alpha \beta}:= U_{\alpha} \cap U_{\beta}$
these bases tranform according to 
$e_{\alpha i|_{U_{\alpha \beta}}} = \sum_{j=1}^{n}(t_{\alpha \beta})_{ij}
e_{\beta j|_{U_{\alpha \beta}}}$. Define a connection $\nabla_{\alpha}$ on
$U_{\alpha}$ by $\nabla_{\alpha}e_{\alpha i}= 0$, for $ 1\leq i \leq n$. Since the 
$t_{\alpha \beta}$ are locally constant, $\nabla_{\alpha} = \nabla_{\beta}$ when restricted to 
$U_{\alpha \beta}$. Therefore there exists a flat connection $\nabla$ defined on $E$ whose 
restriction to $U_{\alpha}$ is $\nabla_{\alpha}$. 

Let $t' \in C^{1}({\cal U'},{\cal G}_{n})$ be another representative of $\theta$,
where ${\cal U}'= \{ U'_{\alpha}: \alpha \in A' \}$. 
It must be shown that $(E,\nabla)$ and $(E',\nabla')$, determined respectively 
by $t$ and $t'$, are isomorphic. There exists a common refinement 
${\cal V}=\{V_{\alpha}: \alpha \in B \}$ 
of ${\cal U}$ and ${\cal U'}$,  
$V_{\alpha} \subset U_{\rho(\alpha)}$ and  $V_{\alpha} \subset U'_{\rho'(\alpha)}$,
such that the respective images $\hat{t}$ and $\hat{t}'$ of 
$t$ and $t'$ in ${\cal C}^{1}({\cal V}, {\cal G}_{n})$ are cohomologous: 
$\hat{t}'_{\alpha \beta} = b_{\alpha}^{-1}\hat{t}_{\alpha \beta}b_{\beta}$,
for some $b\in C^{0}({\cal V},{\cal G}_{n})$.  This merely states that $\hat{t}$
and $\hat{t}'$ define isomorphic vector bundles $\hat{E}$ and $\hat{E}'$, respectively,
over $M$. Recall that the flat connection $\hat{\nabla}$ on $\hat{E}$ is determined by 
requiring the canonical sections $\hat{e}_{\alpha i}$ over $V_{\alpha}$ to be parallel, and 
similarly for $\hat{\nabla}'$ on $\hat{E}'$. The canonical local
sections for $\hat{E}$ and $\hat{E}'$ are related by 
$\hat{e}_{\alpha i} = \sum_{j=1}^{n}(b_{\alpha})_{ij}\hat{e}'_{\alpha j}$. Since the
$b_{\alpha}$ are locally constant, the connections are equal under the canonical
isomorphism of $\hat{E}$ and $\hat{E}'$. That is, $(\hat{E}',\hat{\nabla}') \cong 
(\hat{E},\hat{\nabla})$. 
Furthermore, $E$ is canonically isomorphic to $\hat{E}$. The canonical local sections of these 
two bundles are related by $\hat{e}_{\alpha i} = e_{\rho(\alpha)i|_{V_{\alpha}}}$ and so the
associated flat connections are equal under the canonical isomorphism of the two bundles.
Hence $(E,\nabla) \cong (\hat{E},\hat{\nabla})$ and similarly
$(E',\nabla') \cong (\hat{E}',\hat{\nabla}')$. Therefore, 
$(E',\nabla') \cong (E,\nabla)$ and so $\Psi: H^{1}(M, {\cal G}_{n}) \rightarrow
{\cal E}_{n}$ is well-defined. \\

\hspace{-0.25in}(iii) and (iv) follow directly from the definitions of $\Theta$ and 
$\Psi$ given above. \\
{\bf q.e.d.} 

\begin{corollary} \label{corollary:bijection}
The set of isomorphism classes ${\cal E }_{n}$ is naturally bijective with 
$H^{1}(M, {\cal G}_{n})$.
\end{corollary}

\begin{proposition} \label{proposition:frame}
Let $\nabla$ be a flat connection on the vector bundle $E$ and let $\xi $ denote the
isomorphism class of $(E,\nabla)$. There exists a frame of parallel sections of $E$ 
if and only if $\Theta(\xi) = 1$.
\end{proposition}
{\bf Proof:}\\ 
$\Longrightarrow$ Suppose that $(s_{1},...,s_{n})$ is a frame of parallel sections
of $E$. With respect to any open cover ${\cal U}=\{ U_{\alpha}: \alpha \in A \}$ of $M$,
the transition functions of $s_{\alpha i} := s_{i|_{U_{\alpha}}}$ on non-empty intersections 
$U_{\alpha \beta} :=U_{\alpha} \cap U_{\beta}$
is given by the identity matrix: 
$s_{\alpha i|_{U_{\alpha \beta}}} = \sum_{j=1}^{n}(1_{\alpha \beta})_{ij}
s_{\beta j|_{U_{\alpha \beta}}}$, for $1\leq i\leq n$. 
Therefore $\Theta(\xi) = 1$. \\
$\Longleftarrow$ Suppose that $\Theta(\xi) = 1$. Then there exists an open cover
${\cal U}=\{ U_{\alpha}: \alpha \in A \}$ of $M$ along with local parallel sections $s_{\alpha i}$
on $U_{\alpha}$
such that the transition functions are given by the identity matrix. That is,
on non-empty intersections $U_{\alpha \beta}$, $s_{\alpha i|_{U_{\alpha \beta}}} =
s_{\beta i|_{U_{\alpha \beta}}}$. Therefore there exists a parallel frame $(s_{1},...,s_{n})$
such that $s_{i|_{U_{\alpha}}} = s_{\alpha i}$, $1\leq i \leq n$.\\
{\bf q.e.d.} \\

Note that triviality of the vector bundle $E$ does not imply that
$\Theta(\xi)=1$; $E$ must satisfy the stronger requirement of possessing
a {\it parallel} frame. This is illustrated in the following example. \\

\hspace{-0.3in} {\bf Example 2} \hspace{0.03in} Consider the 
trivial line bundle $E:= S^{1}\times \Re$
over $M=S^{1}$ with the obvious projection. Coordinates on $E$ are given by 
$(\phi,x)$ where $0\leq \phi < 2 \pi$ and $x\in \Re$. $X_{1} := (\phi,1)$ 
defines a frame. Define the flat connection on $E$ by 
\[ \nabla_{\frac{\partial}{\partial \phi}}X_{1} = X_{1} \]
and denote the isomorphism class of $(E,\nabla)$ by $\xi\in {\cal E}_{1}$.
$X(\phi) := e^{-\phi}X_{1}$ is a non-vanishing parallel section but does not extend
over the entire manifold $M$. $E$ does not admit a parallel frame 
and so by Proposition \ref{proposition:frame}, 
$\Theta(\xi) \neq 1$ in $H^{1}(M,{\cal G}_{1})$. 

\begin{corollary} \label{corollary:SC}
$H^{1}(M,{\cal G}_{n})=\{1\}$, for any simply-connected manifold $M$.
\end{corollary}
{\bf Proof:}\\ 
Let $\theta \in H^{1}(M,{\cal G}_{n})$ and choose a representative $(E,\nabla)$
of $\Psi(\theta)$. Since $\nabla$ is a flat connection on $E$ and $M$ is simply-connected,
there exists a frame of parallel sections $s_{i}:M\rightarrow E$ obtained by parallel
transport. By Proposition \ref{proposition:flat} (iii) and Proposition 
\ref{proposition:frame}, $\theta = 1$. \\
{\bf q.e.d.}

\section{The global problem: $rank \hspace{0.03in} \widetilde{W}= n$}

Using the relationship between flat vector bundles and cohomology developed in the 
previous section we now proceed to characterize the regular connections that are globally 
metric. 

\begin{lemma} \label{lemma:posdef}
Let $V$ be a finite-dimensional vector space of symmetric, bilinear forms defined
on some linear space.
If  $V$ contains a positive-definite bilinear form then there exists a basis for 
$V$ comprised of positive-definite bilinear forms.
\end{lemma}
{\bf Proof:}\\ 
Let $e$ be a positive-definite bilinear form in $V$. Extend to a basis
$(e_{1},...,e_{n})$ of $V$, where $e_{1}:=e$. For a sufficiently
small $\epsilon >0$, $b_{i} := e+\epsilon e_{i}$, for $2\leq i \leq n$,
is positive-definite. Set $b_{1}:=e=e_{1}$. Then $(b_{1},...,b_{n})$ is a basis 
of $V$ comprised of positive-definite bilinear forms. \\
{\bf q.e.d.} \\

Consider the natural embedding of groups $\iota: GL(m,\Re) 
\rightarrow GL(n,\Re)$ given by
\[ \iota(A) := \left( \begin{tabular}{c|c} 
$I_{(n-m)\times (n-m)}$ & $0$ \\ \hline
$0$ & $A_{m\times m}$  \end{tabular} \right), \]
for $0\leq m \leq n$. This induces an embedding of sheaves
$\iota:{\cal G}_{m} \rightarrow {\cal G}_{n}$. A 1-cocycle 
$a\in {\cal C}^{1}({\cal U}, {\cal G}_{m})$ is sent to a 1-cocycle $\iota(a) \in
{\cal C}^{1}({\cal U}, {\cal G}_{n})$, which determines an element $[\iota(a)] \in
H^{1}(M,{\cal G}_{n})$. If two cocycles, $a_{1} \in {\cal C}^{1}({\cal U}_{1}, {\cal G}_{m})$ 
and $a_{2} \in {\cal C}^{1}({\cal U}_{2}, {\cal G}_{m})$,
are cohomologous in $H^{1}(M, {\cal G}_{m})$ then $[\iota(a_{1})] = [\iota(a_{2})]$
in $H^{1}(M, {\cal G}_{n})$. Therefore $\iota$ defines a map 
$\hat{\iota}:H^{1}(M, {\cal G}_{m}) \rightarrow H^{1}(M, {\cal G}_{n})$. 
$\hat{\iota}$ is not necessarily injective since cocycles that are not
cohomologous in  $H^{1}(M, {\cal G}_{m})$ may be cohomologous in $H^{1}(M, {\cal G}_{n})$.
Let $H^{1}_{m}$ denote the image of $H^{1}(M, {\cal G}_{m})$ under $\hat{\iota}$. We have
a sequence of subsets
\[ \{1 \} = H^{1}_{0} \subseteq H^{1}_{1} \subseteq \cdots \subseteq H^{1}_{n-1} 
       \subseteq H^{1}_{n} =  H^{1}(M, {\cal G}_{n})  \]
Define the {\it rank} of an element $\sigma \in  H^{1}(M, {\cal G}_{n})$ to be $n-m$,
where $m$ is the least non-negative integer for which $\sigma \in H^{1}_{m}$. 
For instance, $rank \hspace{0.03in} 1 = n$.
This appellation will be justified below.

Now specify $W$ to be the vector bundle of symmetric
elements of $T^{*}M\otimes T^{*}M$ over $M$ and 
let $\nabla$ be a connection on $M$ that is regular with respect to $W$.
Denote the isomorphism class of $(\widetilde{W}, \nabla)$ by $\xi \in {\cal E}_{n}$ and 
set $\tau := \Theta(\xi)$. Let $W_{S}$ be the trivial subbundle of 
$\widetilde{W}$ generated by the global parallel sections of $W$ and define $W_{M}$ to 
be the trivial subbundle of $W_{S}$ generated by the global parallel metrics on $M$.  
The subbundles have the following inclusions:
\[ W_{M} \subseteq W_{S} \subseteq \widetilde{W} \subseteq W \]

\begin{lemma} \label{lemma:dimensions}
If $\nabla$ is a metric connection then $W_{M} = W_{S}$.
\end{lemma}
{\bf Proof:}\\ 
Let $s_{1},...,s_{m}$ be a basis for the space of global parallel sections of $W$  
where $s_{1}$ is a parallel metric. Fix any $x\in M$ and define $V$ to be
the vector space spanned by $\{s_{1}(x),...,s_{m}(x) \}$. $s_{1}(x)$ is a 
positive-definite bilinear form, so by Lemma \ref{lemma:posdef},  there exists a basis 
$(b_{1},...,b_{m})$ of $V$ consisting of positive-definite bilinear forms. It follows that
there exist $m$ linearly independent parallel sections $(g_{1},...,g_{m})$ of $W$ satisfying
$g_{i}(x) = b_{i}$, for $1\leq i \leq m$. Since the $b_{i}$ are positive-definite, the $g_{i}$
are parallel metrics. Therefore $rank \hspace{0.03in} W_{M} \geq  m =
rank \hspace{0.03in} W_{S}$. But $W_{M}$ is a subbundle of $W_{S}$ and so 
$W_{M} = W_{S}$. \\
{\bf q.e.d.}

\begin{lemma} \label{lemma:parallelmetric} 
$rank \hspace{0.03in} \tau \geq rank \hspace{0.03in} W_{M} $. 
\end{lemma}
{\bf Proof:}\\ 
If $rank \hspace{0.03in} W_{M} =0$ then the lemma holds trivially. Therefore 
suppose $rank \hspace{0.03in} W_{M} = m > 0$; there exist $m$ linearly 
independent parallel metrics $s_{i}$ on $M$, $1\leq i \leq m$. Let $s$ be any parallel metric
on $M$, say $s=s_{1}$. Define a smooth inner product $<,>$ on the fibres of 
$\widetilde{W}$ by
\[ <h,h'> := \sum_{\mu, \nu, \kappa, \lambda } h_{\mu \nu}h'_{\kappa \lambda}
s^{\mu \kappa}s^{\nu \lambda} \]
for $h,h' \in \widetilde{W}_{x}$ and $x\in M$. Let $\widetilde{W}_{1}$ be the rank
$m$ subbundle of $\widetilde{W}$ generated by the $s_{i}$ and let 
$\widetilde{W}_{2}$ be the subbundle perpendicular to $\widetilde{W}_{1}$ with
respect to the inner product. $\widetilde{W}$ is the direct sum of the 
two subbundles: $\widetilde{W}=\widetilde{W}_{1} \oplus \widetilde{W}_{2}$.

We claim that $\widetilde{W}_{1}$ and $\widetilde{W}_{2}$ are invariant with respect to 
parallel transport. This is clear for $\widetilde{W}_{1}$ since it is generated by global
parallel metrics. Consider the parallel transport $h$ 
along a curve $\gamma:[0,1] \rightarrow M$ of an element $h(\gamma(0)) \in \widetilde{W}_{2}$.
At $\gamma(0)$, 
\[ <s_{i},h> := \sum_{\mu, \nu, \kappa, \lambda} (s_{i})_{\mu \nu}h_{\kappa \lambda}
s^{\mu \kappa}s^{\nu \lambda} \]
equals zero, by definition of $\widetilde{W}_{2}$, for each $1\leq i \leq m$. The 
covariant derivative of the right hand side of the equation with respect to 
$\gamma_{*}(\partial / \partial u)$ at $\gamma (u)$ vanishes because each of the four
contracted tensors is parallel along $\gamma$. Therefore the left hand side is a constant
along $\gamma$ and must be identically equal to zero. This means that $h(\gamma(1))$
is an element of $\widetilde{W}_{2}$, establishing the claim.

For $y\in M$ let $(e_{y1},...,e_{yn})$ be a basis of $\widetilde{W}_{y}$ where
$e_{yi}=s_{i}(y) \in \widetilde{W}_{1}$ for $1\leq i \leq m$ and $e_{yi} \in 
\widetilde{W}_{2}$ for $m+1 \leq i \leq n$. There exists an open neighbourhood $U_{y}$ of $y$ 
and a basis of parallel sections $h_{yi}$ of $\tilde{\pi}^{-1}(U_{y})$ such that
$h_{yi}(y) = e_{yi}$ for $1\leq i \leq n$. It follows that $h_{yi} = s_{i|_{U_{y}}}$
for $1\leq i \leq m$, and since $\widetilde{W}_{2}$ is invariant under parallel transport,
$h_{yi}:U_{y} \rightarrow \widetilde{W}_{2}$ for
$m+1 \leq i \leq n$. Consider two open neighbourhoods $U_{y}$ and $U_{z}$ in ${\cal U} := 
\{ U_{x} : x\in M \}$ with non-empty intersection $U_{yz} := U_{y}\cap U_{z}$.
The element $t\in {\cal C}^{1}({\cal U},{\cal G}_{n})$ describing the change of basis 
$h_{yi|_{U_{yz}}} = \sum_{j=1}^{n}(t_{yz})_{ij}h_{zj|_{U_{yz}}}$ has the form
$t=\iota(a)$ for some $a\in {\cal C}^{1}({\cal U},{\cal G}_{n-m})$:
\[ t_{yz} = \left( \begin{tabular}{c|c} 
$1_{yz}$ & $0$ \\ \hline
$0$ & $a_{yz}$  \end{tabular} \right) \]
Therefore $\tau = [t] \in H^{1}_{n-m}$ and so 
$rank \hspace{0.03in} \tau \geq m =rank \hspace{0.03in} W_{M} $. \\
{\bf q.e.d.} \\

Let $rank \hspace{0.03in} \tau = m$. 
Then we can find an open cover ${\cal U}= \{U_{\alpha} : \alpha \in A \}$
of $M$ and associated basis $(h_{\alpha 1}, ..., h_{\alpha n})$ of parallel sections
on $U_{\alpha}$ such that the transformation of the bases over the 
non-empty intersection $U_{\alpha \beta} := U_{\alpha} \cap U_{\beta}$ of two open sets
in ${\cal U}$, 
\[ h_{\alpha i|_{U_{\alpha \beta}}} = \sum_{j=1}^{n}(t_{\alpha \beta})_{ij}
   h_{\beta j|_{U_{\alpha \beta}}}  \]
is given by $t=\iota(a)$ for some $a\in {\cal C}^{1}({\cal U},{\cal G}_{n-m})$. 
Therefore, restricted to $U_{\alpha \beta}$, $h_{\alpha i|_{U_{\alpha \beta}}}
= h_{\beta i|_{U_{\alpha \beta}}}$ for all $1\leq i \leq m$.
It follows that there exist $m$ independent global parallel sections $s_{1},...,s_{m}$ 
whose restriction to $U_{\alpha}$ is $s_{i|_{U_{\alpha}}} = h_{\alpha i }$. 
This demonstrates the following lemma.

\begin{lemma} \label{lemma:wsandtau}
$rank \hspace{0.03in} W_{S} \geq rank \hspace{0.03in} \tau$.
\end{lemma}

The theorem below relates the number of independent parallel metrics of a regular 
metric connection $\nabla$ to the associated cohomology. 
 
\begin{theorem} \label{theorem:dimension}
If $\nabla$ is a regular metric connection then 
$rank \hspace{0.03in} W_{M} = rank \hspace{0.03in} W_{S} = rank \hspace{0.03in} \tau $.
\end{theorem} 
{\bf Proof:}\\ 
From Lemmas \ref{lemma:dimensions}, \ref{lemma:parallelmetric} and \ref{lemma:wsandtau},
$rank \hspace{0.03in} W_{M} = rank \hspace{0.03in} W_{S} \geq rank \hspace{0.03in} \tau
\geq rank \hspace{0.03in} W_{M}$.  \\
{\bf q.e.d.} \\

Let $W_{\nabla}$ denote the trivial subbundle of $W_{S}$ generated by the $h_{\alpha i }$
for $1\leq i \leq m$ and $\alpha \in A$.  

\begin{proposition} \label{proposition:welldefined}
$W_{\nabla}$ is well-defined.
\end{proposition}
{\bf Proof:}\\ 
Suppose that $\tau $ is represented by another element $t'=\iota(a')$, where $a' \in
{\cal C}^{1}({\cal U'},{\cal G}_{n-m})$, with respect to an open cover ${\cal U}'$ of $M$.
Let $s'_{1},...,s'_{m}$ be the associated independent global parallel sections. 
We must show that the bundle generated by the $s'_{i}$ equals the
bundle generated by the $s_{i}$. Assume otherwise. Then there is a global parallel section
$s_{0}:=s'_{l}$ which is not in the span of the $s_{i}$. Consider an open set $U_{\alpha}$
in ${\cal U}$. Recall that $\{ s_{1|_{U_{\alpha}}},...,
s_{m|_{U_{\alpha}}},h_{\alpha m+1},...,h_{\alpha n} \}$ is a basis for the space of parallel
sections on $U_{\alpha}$. However, the set of parallel sections 
$\{ s_{0|_{U_{\alpha}}},s_{1|_{U_{\alpha}}},..., s_{m|_{U_{\alpha}}},
h_{\alpha m+1},...,h_{\alpha n} \}$ on $U_{\alpha}$ is linearly dependent,
consisting of $n+1$ elements. 
Moreover, for some $k \in \{m+1,..., n\}$, $h_{\alpha k}$ can be written as
a linear combination of the other elements in the set:
\[ h_{\alpha k} = \sum_{i=0}^{m}c_{i}s_{i|_{U_{\alpha}}}+
                  \sum_{i=m+1}^{k-1}c_{i}h_{\alpha i} + 
                  \sum_{i=k+1}^{n}c_{i}h_{\alpha i}     \]
Therefore
\[ (\bar{h}_{\alpha 1} ,..., \bar{h}_{\alpha n}) := 
 (s_{0|_{U_{\alpha}}},...,s_{m|_{U_{\alpha}}}, h_{\alpha m+1},...,h_{\alpha k-1},
 h_{\alpha k+1},..., h_{\alpha n}) \]
is a basis of local parallel sections on $U_{\alpha}$.   
The transformation of these bases over the 
non-empty intersection $U_{\alpha \beta} := U_{\alpha} \cap U_{\beta}$ of two open sets
in ${\cal U}$, 
\[ \bar{h}_{\alpha i|_{U_{\alpha \beta}}} = \sum_{j=1}^{n}(\bar{t}_{\alpha \beta})_{ij}
   \bar{h}_{\beta j|_{U_{\alpha \beta}}}  \]
is given by $\bar{t}=\iota(\bar{a})$ for some 
$\bar{a} \in {\cal C}^{1}({\cal U},{\cal G}_{n-m-1})$.
This contradicts the stated rank of $\tau$. \\
{\bf q.e.d.} 

\begin{proposition} \label{proposition:equality}
If $\nabla$ is a metric connection then $W_{\nabla} = W_{S} = W_{M}$.
\end{proposition}
{\bf Proof:}\\ 
Suppose $\nabla$ is a metric connection. By Theorem \ref{theorem:dimension}, 
$rank \hspace{0.03in} W_{S} = m$. But $W_{\nabla}$ is a rank $m$ subbundle of $W_{S}$
and so $W_{\nabla} = W_{S}$, which equals $W_{M}$ by Lemma \ref{lemma:dimensions}.\\
{\bf q.e.d.} \\

We may now determine when a regular connection is a metric connection.

\begin{theorem} \label{theorem:main}
Let $\nabla$ be a regular connection on $M$.
Then  $\nabla$ is a metric connection if and only if  $W_{\nabla}$ contains a 
positive-definite bilinear form. 
\end{theorem}
{\bf Proof:}\\ 
$\Longrightarrow$ Suppose that $\nabla$ is a metric connection. 
By Proposition \ref{proposition:equality}, $W_{\nabla} = W_{M}$.
Let $g$ be any parallel metric. Then for any $x\in M$, 
$g(x) \in W_{M}$ is a positive-definite bilinear form in $W_{\nabla}$. \\
$\Longleftarrow$ Suppose that $W_{\nabla}$ contains a positive-definite bilinear form $e$,
lying in the fibre of $W_{\nabla}$ over $x\in M$, say. Then there exist constants $c_{i}$
such that $e=c_{1}s_{1}(x)+\cdots+c_{m}s_{m}(x)$. $s:=c_{1}s_{1}+\cdots+c_{m}s_{m}$
is a parallel metric. \\
{\bf q.e.d.} \\

It is worthwhile comparing Theorem \ref{theorem:main} and Corollary \ref{corollary:local}.
$W_{\nabla}$ plays the analogous role, in the global context, to 
$\widetilde{W}_{x}$ in the local problem. In both cases, metricity of the
connection is determined by the existence of a positive-definite bilinear form in the relevent 
space. There is a fundamental distinction, however:
$\widetilde{W}_{x}$ is obtained by algebraic solution to homogeneous linear systems, 
whereas $W_{\nabla}$, on the other hand, must be found by integration of 
differential equations. \\

\hspace{-0.3in} {\bf Example 3} Consider the symmetric 
connection $\nabla$ on the punctured plane $M=\Re^{2}-\{(0,0)\}$, given by 
$\Gamma^{r}_{\theta \theta}= -k^{2}r$, 
$\Gamma^{\theta}_{r \theta} =  \Gamma^{\theta}_{\theta r} = \frac{1}{r}$ and 
all other Christoffel symbols equal to zero.
$0<r<\infty$ is the radial distance, $0\leq \theta<2\pi$ is the planar angle 
and $k$ is a non-half-integer constant.
$\nabla$ is flat and so $\widetilde{W} = W^{(0)}$  is a rank three vector bundle over $M$
spanned by 
\[ \begin{array}{lll}
   X_{1} & = & dr \otimes dr \\
   X_{2} & = & d\theta \otimes d\theta \\
   X_{3} & = & dr \otimes d\theta + d\theta \otimes dr
   \end{array} \]
A section $h=f_{1}X_{1}+f_{2}X_{2}+f_{3}X_{3}$ of $\widetilde{W}$ is parallel if
$\nabla h = 0$, which is equivalent to a system of six partial 
differential equations in the three unknown functions $f_{1},f_{2}$, and $f_{3}$. These are
readily integrated to yield three independent local parallel sections:
\[ \begin{array}{lll}
   h_{1} & = & X_{1}+k^{2}r^{2}X_{2} \\
   h_{2} & = & k^{-1}sin(2k\theta)X_{1} - kr^{2}sin(2k\theta) X_{2} +
               rcos(2k\theta)X_{3} \\
   h_{3} & = & k^{-1}cos(2k\theta)X_{1} - kr^{2}cos(2k\theta) X_{2} -
               rsin(2k\theta)X_{3} 
   \end{array} \]
No non-zero linear combination of $h_{2}$ and $h_{3}$ can be extended over all of $M$.
Therefore $\nabla$ has exactly one global parallel metric $h_{1}$, up to constant multiples.
$W_{\nabla} = W_{S} = W_{M}$ is generated by $h_{1}$ and 
$rank \hspace{0.03in} \tau = rank \hspace{0.03in} W_{M}=1 $. 

From the perspective of cohomology, let $U_{1}:=\Re^{2}-\{ (x,0): x\geq 0 \}$ and 
$U_{2}:=\Re^{2}-\{ (x,0): x\leq 0 \}$ define an open cover ${\cal U}$ of $M$.
$U_{1,2} := U_{1}\cap U_{2} = H^{+}\cup H^{-}$ where $H^{+}$ (resp. $H^{-}$) is the upper 
(resp. lower) half plane.  Define parallel sections on $U_{1}$ and $U_{2}$ by 
\[ \begin{array}{lllll}
   h_{1i}(r,\theta)   & := & h_{i|_{U_{1}}} (r,\theta) & 
                             \hspace{0.5in} &  0 <\theta < 2\pi\\
                         
   h_{2i}(r,\theta')   & := & h_{i|_{U_{2}}} (r,\theta') & 
                             \hspace{0.5in} &  \pi < \theta' < 3\pi
   \end{array} \]
for $1\leq i \leq 3$, where $0< \theta < 2\pi$ and $\pi<\theta' <3\pi$ are coordinates on
$U_{1}$ and $U_{2}$, respectively. The transition between the $h_{1i}$ and  $h_{2j}$ 
on $U_{1,2}$ is given by 
$h_{1i|_{U_{1,2}}} = \sum_{j=1}^{3}(t_{12})_{ij}h_{2j|_{U_{1,2}}}$ for
 $t\in C^{1}({\cal U},{\cal G}_{3})$ defined by
 \[ \begin{array}{llc}
   t_{12|_{H^{-}}}    & := & I_{3\times 3} \\
                      &    &  \\ 
   t_{12|_{H^{+}}}    & := &  \left( \begin{array}{ccc}
                             1 \hspace{0.2in} & 0 & 0 \\
                             0 \hspace{0.2in} & cos(4k\pi) & -sin(4k\pi) \\
                             0 \hspace{0.2in} & sin(4k\pi) & cos(4k\pi)
                             \end{array}
						      \right)
   \end{array} \]
The image of $t$ in $H^{1}(M,{\cal G}_{3})$ defines the cohomology class associated
to the connection: $\tau = [t] \in H^{1}_{2} \subseteq H^{1}(M,{\cal G}_{3})$.  \\

An example of a flat non-metric connection on the torus may be found in 
\cite{ff}, Example 4.2, pg. 211.

\section{Analytic versus smooth connections}

In this final section we consider the differences between smooth and analytic connections.

\begin{proposition} \label{proposition:smooth}
Let $\nabla$ be a regular, locally metric connection on a simply-connected manifold $M$. 
Then $\nabla$ is a metric connection. Furthermore, there exist 
$rank \hspace{0.03in} \widetilde{W}$
independent parallel metrics on $M$.
\end{proposition}
{\bf Proof:}\\
By Corollary \ref{corollary:SC}, $\tau = 1$ and so $W_{\nabla} = W_{S} = \widetilde{W}$,
which contains a positive-definite bilinear form since $\nabla$ is a locally metric
connection. By Theorem \ref{theorem:main}, $\nabla$ is metric and by Lemma 
\ref{lemma:dimensions}, 
$rank \hspace{0.03in} W_{M} = rank \hspace{0.03in} \widetilde{W}$.\\
{\bf q.e.d.} \\

The regularity requirement is superfluous, however, for analytic connections.

\begin{lemma} \label{lemma:analytic-regular}
If $\nabla$ is an analytic connection on an analytic manifold $M$
then $\widetilde{W}$ is a vector bundle over $M$.
\end{lemma}
{\bf Proof:}\\ 
For every point $x \in M$ there exists an open neighbourhood $U_{x}$ of $x$ and 
a (possibly empty) basis of local analytic parallel metrics $h_{xi}$ , $1\leq i \leq N_{x}$,
on $U_{x}$, which has maximal cardinality; if $V$ is any other open neighbourhood of $x$ and 
$g_{i}$, $1\leq i\leq N$ a basis of local analytic parallel metrics on $V$, then $N \leq N_{x}$.
Without loss of generality we may assume that the $U_{x}$ are diffeomorphic to 
$\Re^{n}$
and that the intersection of any two such open neighbourhoods is also diffeomorphic
to $\Re^{n}$, since a good refinement of ${\cal U} := \{ U_{x}: x\in M\}$ exists.

Consider two open sets $U_{x}$ and $U_{y}$ in ${\cal U}$ with non-empty intersection
$U_{x}\cap U_{y}$. Each $h_{xi}$ extends, by analytic continuation, to an
analytic metric $h'_{xi}$ on $U_{x}\cup U_{y}$. Since $\nabla$ is analytic, so is
$\nabla h'_{xi}$. Restricted to $U_{x}$, 
$\nabla h'_{xi|_{U_{x}}} = 0$ and hence, by analyticity, 
$\nabla h'_{xi} = 0$ on all of $U_{x}\cup U_{y}$. Therefore $h'_{xi|_{U_{y}}}$ are 
$N_{x}$ independent local parallel
metrics on $U_{y}$. This shows that $N_{x} \leq N_{y}$ by the maximal cardinality
requirement. By symmetry $N_{y} \leq N_{x}$ and so $N_{x}= N_{y}$. Since $M$ is assumed to
be a connected manifold $N_{x}$ is independent of $x$ and equals the dimension of the 
fibres of $\widetilde{W}$ over $M$. 
The fibre bundle structure of $\widetilde{W}$ is supplied by the transition functions between
the local bases $h_{xi}$.\\
{\bf q.e.d.} 

\begin{corollary} \label{corollary:analytic}
Let $\nabla$ be an analytic locally metric connection on an analytic
simply-connected manifold $M$. Then $\nabla$ is a metric connection. Furthermore, there exist
$rank \hspace{0.03in} \widetilde{W}$ independent parallel metrics on $M$.
\end{corollary}

Analyticity is essential in the above corollary; it is not true for smooth 
connections as the following example shows. \\

\hspace{-0.3in} {\bf Example 4} 
Let
 \[ f(x) := \left\{ 
   \begin{array}{lll}
   a+ e^{\frac{-1}{(x-x_{0})^{2}}} & \hspace{0.5in} & x < x_{0} \\
   a  & \hspace{0.5in} &  x \geq x_{0} 
   \end{array} \right. \]
and 
\[ h(x) := \left\{ 
   \begin{array}{lll}
   b & \hspace{0.5in} & x < x_{1} \\
   b+ e^{\frac{-1}{(x-x_{1})^{2}}}  & \hspace{0.5in} &  x \geq x_{1} 
   \end{array} \right. \]
where $a,b,x_{0}$ and $x_{1}$ are constants satisfying $a,b >0$, $a\neq b$ and $x_{0} < x_{1}$.
Define the smooth connection
$\nabla$ on $\Re^{2}$ by 
\[ \begin{array}{llll}
\Gamma^{x}_{yy} &    :=  &    \left\{
   \begin{array}{c}
   -\frac{f'}{2}  \\
   0   \\
   -\frac{h'}{2} 
   \end{array} \right.  & \hspace{0.4in}
  \begin{array}{l}
   x< x_{0} \\
   x_{0}\leq x \leq x_{1}   \\
   x > x_{1}
   \end{array}       \\
 \hspace{.001in} & \hspace{.001in} & \hspace{.001in} & \hspace{.001in}  \\  
\Gamma^{y}_{xy} = \Gamma^{y}_{yx} &    :=  &     \left\{ 
   \begin{array}{c}
   \frac{f'}{2f}  \\
   0   \\
   \frac{h'}{2h} 
   \end{array} \right.  &  \hspace{0.4in}
  \begin{array}{l}
   x< x_{0} \\
   x_{0}\leq x \leq x_{1}   \\
   x> x_{1}
   \end{array}   
\end{array}  \]
and all other Christoffel symbols equal to zero.
$\nabla$ is a locally metric connection: for $x< x_{1}$, the parallel metrics are
$c(dx^{2} + f(x)dy^{2})$,
and for $x > x_{0}$, the parallel metrics are 
$c(dx^{2} + h(x)dy^{2})$, where $c$ is an arbitrary positive constant. 
Since $a\neq b$, elements of these two sets cannot be joined together to give a global 
parallel metric. Therefore $\nabla$ is not a metric connection  on $\Re^{2}$. 

Observe that $\widetilde{W}$ is not a fibre bundle as the dimension of the
fibres of $\widetilde{W}$ varies over the manifold: on 
$x \leq x_{0}$,  $dim\hspace{0.03in} \widetilde{W}_{x} = 1$,
on $x_{0} < x < x_{1}$, $dim\hspace{0.03in} \widetilde{W}_{x} = 3$ and on $x \geq x_{1}$, 
$dim\hspace{0.03in} \widetilde{W}_{x} = 1$.
This behaviour prevents the formulation of the global existence problem 
in terms of a sheaf cohomology of groups.

\end{document}